\documentclass[12pt]{article}

\usepackage[english]{babel} 
\usepackage{times}
\usepackage{graphicx}
\usepackage{amsmath}
\usepackage{amsfonts}
\usepackage{amssymb}
\usepackage{amsthm}
\usepackage{array}
\usepackage{verbatim}
\usepackage{epsfig}
\usepackage[dvipsnames]{xcolor}
\usepackage{latexsym}
\usepackage{float,booktabs}
\usepackage{dcolumn}
\usepackage{bm}
\usepackage[list=true]{subcaption}
\usepackage{tikz}
\usepackage{textpos}
\usepackage{helvet}
\usepackage{array}
\usepackage{siunitx}
\usepackage{framed}
\usepackage{mathtools}
\usepackage{multicol}
\usepackage{enumitem}
\usepackage{dsfont}
\usepackage{calc}
\usepackage[ruled,vlined]{algorithm2e}

\def\norm#1{\|#1\|}

\title{Estimation of distributions via multilevel Monte Carlo with stratified sampling}

\author{S{\o}ren Taverniers\footnote{Department of Energy Resources Engineering, Stanford University, 367 Panama St, Stanford, CA 94305, USA, \texttt{tartakovsky@stanford.edu}} \and Daniel M. Tartakovsky\footnotemark[\value{footnote}]}

\begin{document}
\maketitle

\begin{abstract}

We design and implement a novel algorithm for computing a multilevel Monte Carlo (MLMC) estimator of the cumulative distribution function of a quantity of interest in problems with random input parameters or initial conditions. Our approach combines a standard MLMC method with stratified sampling by replacing standard Monte Carlo at each level with stratified Monte Carlo with proportional allocation. We show that the resulting stratified MLMC algorithm is  more efficient than its standard MLMC counterpart, due to the reduction in variance at each level provided by the stratification of the random parameter's domain. A smoothing approximation for the indicator function based on kernel density estimation yields a more efficient algorithm compared to the typically used polynomial smoothing. The difference in computational cost between the smoothing methods depends on the required error tolerance.

\end{abstract}

\section{Introduction and Motivation}
\label{sec:intro}

Simulation of many complex systems, such as subsurface flows in porous media~\cite{Hill2007,Tartakovsky2013} or reaction initiation in heterogeneous explosives~\cite{Sen2018}, is complicated by a lack of information about key properties such as permeability or initial porosity. Uncertainty in the medium's properties or initial state propagates into uncertainty in predicted quantities of interest (QoIs), such as mass flow rate or the material's temperature. 

Probabilistic methods, which treat an uncertain input or initial state of the system as a random variable, provide a natural venue to quantify predictive uncertainty in a QoI. These techniques render the QoI random as well, i.e., it takes on values that are distributed according to some probability density function (PDF). Such approaches include stochastic finite element methods (FEMs), which characterize the random parameter fields in terms of a finite set of random variables, e.g., via a spectral representation or a Karhunen-Lo\`eve expansion. This finite set of random variables defines a finite-dimensional outcome space on which the solution to the resulting stochastic partial differential equation (PDE) is defined. Examples of stochastic FEMs include stochastic Galerkin, which expands the solution of a stochastic PDE in terms of orthogonal basis functions, and stochastic collocation, which samples the random parameters at predetermined values or ``nodes''~\cite{Xiu2010}. While such methods perform well when the number of stochastic parameters (aka ``stochastic dimension'') is low and these parameters exhibit a long correlation length, for many stochastic degrees of freedom and short correlation lengths their performance decreases dramatically~\cite{Taverniers2017}. In addition, since nonlinearity degrades the solution's regularity in the outcome space, stochastic FEMs also struggle with solving highly nonlinear problems~\cite{barajas-2016-stochastic}. Another class of probabilistic techniques involves the derivation of deterministic equations for the statistical moments~\cite{Jarman2003,Winter2003} or PDF~\cite{Lichtner2003,Tartakovsky2009} of the QoI. While these methods do not suffer from the ``curse of dimensionality'', they require a closure approximation for the derived moment or PDF equations. Such closures often involve perturbation expansions of relevant quantities into series in the powers of the random parameters' variances, which limits the applicability of such methods to parameters with low coefficients of variation.

Monte Carlo simulations (MC)~\cite{Fishman1996} remain the most robust and straightforward way to solve PDEs with random parameters/initial conditions. The method samples the random variables from their distribution, solves the deterministic PDE for each realization, and computes the resulting statistics of the QoI. While combining a nonintrusive character with a convergence that is independent of the stochastic dimension, MC converges very slowly: the standard deviation of the MC estimator for the QoI's expectation value is inversely proportional to $\sqrt{N}$ where $N$ is the number of realizations~\cite{Fishman1996}. While this drawback spurred the development of alternative probabilistic methods such as those listed above, efforts to combine MC with the multigrid concept by Heinrich~\cite{Heinrich1998,Heinrich2001} and later by Giles~\cite{Giles2008} have sparked renewed interest in MC under the form of the multilevel Monte Carlo (MLMC) method. MLMC aims to achieve the same solution error as MC but at a lower computational cost by correcting realizations on a coarse spatial grid with sampling at finer levels of discretization. By sampling predominantly at the coarsest levels where samples are cheaper to compute, it aims to outperform MC which only performs realizations on the finest spatial grid. The related technique of \emph{multifidelity} Monte Carlo (also referred to as \emph{solver}-based MLMC as a opposed to traditional \emph{grid}-based MLMC) generalizes this approach by using models of varying fidelities and speeds on the different levels~\cite{Muller2014, Peherstorfer2016, OMalley2018}.

Most work on MLMC has been aimed at computing estimators for expectation values and variances of QoIs (see, e.g.,~\cite{Cliffe2011,Mishra2012,Muller2013,Lu2018}). However, the information contained in these first two moments is insufficient to understand, e.g., the probability of rare events. This task requires knowledge of the QoI's full PDF, or equivalently, its cumulative distribution function (CDF)~\cite{Tversky1992,Goovaerts2001,Boso2016}. Assuming the QoI $Q$ to be a function of a continuous random input variable $W$, i.e., $Q=Q(W)$, where $W:\Omega\rightarrow \mathbb{R}$ is a measurable function with $\Omega$ the sample space, the CDF $F$ of $Q$ at a point $q\in \mathbb{R}$ is given by
\begin{align}
  F(q) = \mathbb{E}[\mathcal{I}_{(-\infty,q]}(Q)] = \int_{-\infty}^{\infty} \mathcal{I}_{(-\infty,q]}(s)\ f(s)\ \mathrm d s = \int_{-\infty}^{q} f(s)\ \mathrm d s,
\end{align}
where the indicator function $\mathcal{I}_{(-\infty,q]}(s)$ is defined as
\begin{align}
   \mathcal{I}_{(-\infty,q]}(s) = 
   \begin{cases}   
       1 & \text{for }s\in(-\infty,q] \\
       0 & \text{for }s\in(q,+\infty).
   \end{cases}
\end{align}
Application of MLMC to the estimation of distributions only started a few years ago. Giles et al.~\cite{Giles2015} developed an algorithm to estimate PDFs and CDFs using the indicator function approach, while Bierig et al.~\cite{Bierig2016} approximated PDFs via a truncated moment sequence and the method of maximum entropy. Both approaches were used to approximate CDFs of molecular species in biochemical reaction networks~\cite{Wilson2016}. Elfverson et al.~\cite{Elfverson2016} used MLMC to estimate failure probabilities, which are single-point evaluations of the CDF. Lu et al.~\cite{Lu2016} obtained a more efficient MLMC algorithm by calibrating the polynomial smoothing of the indicator function in~\cite{Giles2015} to optimize the smoothing bandwidth for a given value of the error tolerance, thereby achieving a faster variance decay with level, and by enabling an a posteriori switch to MC if the latter turned out to have a lower computational cost. Krumscheid et al.~\cite{Krumscheid2018} developed an algorithm for approximating general parametric expectations, including CDFs and characteristic functions, and simultaneously deriving robustness indicators such as quantiles and conditional values-at-risk.

To reduce the computational cost of MLMC further, the multigrid approach can be combined with a sampling strategy at each level that is more efficient than standard MC. For example, quasi-Monte Carlo uses quasi-random, rather than random or pseudo-random, sequences to achieve faster convergence than MC, and has been used to speed up the MLMC computation of the mean system state~\cite{Kuo2017}. Furthermore, a number of so-called ``variance reduction'' techniques have been developed to obtain estimators with a lower variance than MC for the same number of realizations. These include stratification, antithetic sampling, importance sampling, and control variates~\cite{Fishman1996}. Recently, importance sampling was incorporated into MLMC for a more efficient computation of expectation values of QoIs~\cite{Kebaier2018}. Ullman et al.~\cite{Ullmann2015} estimated failure probabilities for rare events by combining subset simulation using Markov chain Monte Carlo~\cite{Au2001} with multilevel failure domains defined on a hierarchy of discrete spatial grids with decreasing mesh sizes. We propose to divide the domain of a random input parameter or initial condition using stratified Monte Carlo to achieve variance reduction at each discretization level, and to estimate the CDF of an output QoI via the resulting ``stratified'' Multilevel Monte Carlo approach.

Section~\ref{sec:MLMC} provides a mathematical formulation of MLMC and sMLMC estimators of CDFs and the cost and error associated with computing them. Section~\ref{sec:Results} discusses two testbed problems for assessing the performance of MLMC and sMLMC compared to standard MC. Conclusions and future research directions are reserved for Section~\ref{sec:Concl}.

\section{Multilevel Monte Carlo for CDFs}
\label{sec:MLMC}

Usually the QoI $Q$ cannot be simulated directly. Instead, one discretizes $Q$ on a spatial grid $\mathcal{T}_M$ and introduces a sequence of random variables, $Q_M$ where $M$ is the number of cells in $\mathcal{T}_M$, which converges to $Q$ as $M$ increases. We assume that $Q_M$ converges both in the mean and in the sense of distribution to $Q$ as $M\rightarrow\infty$, i.e., 
\begin{align}
  \mathbb{E}[Q_M - Q] = \mathcal{O}(M^{-\alpha_1}), \quad \mathbb{E}[\mathcal{I}_{(-\infty,q]}(Q_M) - \mathcal{I}_{(-\infty,q]}(Q)] = \mathcal{O}(M^{-\alpha_2}),
\end{align}
as $M\rightarrow \infty$ for $\alpha_1,\alpha_2\in\mathbb{R}$ independent of $M$ and $q$. We approximate the statistics of $Q_M$ rather than of $Q$, i.e., our goal is to estimate the CDF $F_M$ of $Q_M$ on some compact interval $[a,b]\subset\mathbb{R}$. At each point $q\in [a,b]$, $F_M(q)$ is given by
\begin{align}
  F_M(q) = \mathbb{E}[\mathcal{I}_{(-\infty,q]}(Q_M)].
\end{align}
We assume the PDF $f_M$ of $Q_M$ to be at least $d$-times continuously differentiable on $[a-\xi_0,b+\xi_0]$ for some $d\in\mathbb{N}_0$ and $\xi_0>0$. Then, given the values of $F_M$ at a set of $S+1$ equidistant points $\mathcal{S}_h=\{a=q_0 < q_1 < \dots < q_S=b\}$ with separation distance $h$, we  estimate $F_M(q)$ at any point $q\in [a,b]$ via a piecewise polynomial interpolation of degree $\max(d,1)$. The resulting approximation $F_{h,M}(q)$ of $F_M(q)$ is given by
\begin{align}\label{Fhm}
   F_{h,M}(q) = \sum_{n=0}^S \mathbb{E}[\mathcal{I}_{(-\infty,q_n]}(Q_M)] \phi_n(q) \equiv \sum_{n=0}^S \mathbb{E}[\mathcal{I}_n(Q_M)]\ \phi_n(q),
\end{align}
where $\phi_n$ ($n=0,\dots,S$) are, e.g., Lagrange basis polynomials or, in the case of cubic spline interpolation, third-degree polynomials. The goal is to find an unbiased estimator $\hat{\mathcal{I}}_{n,M}$ of $\mathbb{E}[\mathcal{I}_n(Q_M)]$, so that an estimator $\hat{F}_{h,M}(q)$ of $F_{h,M}(q)$ is computed as
\begin{align}
   \hat{F}_{h,M}(q) = \sum_{n=0}^S \hat{\mathcal{I}}_{n,M} \phi_n(q).
\end{align}
The estimator $\hat{F}_{h,M}$ converges to the CDF $F$ of the original QoI $Q$ as $M\rightarrow \infty$ (to reduce the discretization error) and its variance decreases (to reduce the sampling error).

\subsection{Standard Monte Carlo (MC)}
\label{subsec:MC}

The MC estimator for $\mathbb{E}[\mathcal{I}_n(Q_M)]$ based on $N_\text{MC}$ independent samples of $Q_M$ is defined by
\begin{align}
  \hat{\mathcal{I}}_{n,M}^\text{MC} = \frac{1}{N_\text{MC}}\sum_{j=1}^{N_\text{MC}} \mathcal{I}_n(Q_M^{(j)}),
\end{align}    
where $Q_M^{(j)}$ is the $j^{\text{th}}$ sample of $Q_M$. Let $\hat{F}_{h,M}^\text{MC}$ denote the MC estimator of $F_{h,M}$, and $F_h(q) = \sum_{n=0}^S \mathbb{E}[\mathcal{I}_n(Q)]\ \phi_n$ be a piecewise polynomial interpolation of $F(q)$ given its value at a set of points $\{q_n\}=\{a=q_0 < q_1 < \dots < q_S=b\}$ with $n=0,\dots,S$. The error between the two, $e(\hat{F}_{h,M}^\text{MC})\equiv \hat{F}_{h,M}^\text{MC} - F_h$, consists of a sampling part related to estimating $\mathbb{E}[\mathcal{I}_n(Q)]$ and a discretization part related to approximating $Q$ by $Q_M$. Specifically, the mean squared error (MSE) of $\hat{F}_{h,M}^\text{MC}$ is bounded by
\begin{align}\label{MSE_MC}
   \mathbb{E}[\norm{e(\hat{F}_{h,M}^\text{MC})}_{\infty}^2] &\leq \underbrace{\mathbb{E}[\norm{\hat{F}_{h,M}^\text{MC} - \mathbb{E}[\hat{F}_{h,M}^\text{MC}]}_{\infty}^2]}_{e_1^\text{MC}} + \underbrace{\norm{F_{h,M} - F_h}^2}_{e_2^\text{MC}} \\\nonumber &\leq \max_{0\leq n\leq S} N_\text{MC}^{-1} \mathbb{V}[\mathcal{I}_n(Q_M)] + \max_{0\leq n\leq S} |\mathbb{E}[\mathcal{I}_n(Q_M) - \mathcal{I}_n(Q)]|^2,
\end{align}
where $\norm{\cdot}_{\infty}$ denotes the $L^{\infty}$ norm and $\mathbb{V}[\cdot]$ refers to the variance operator. Here $e_1^\text{MC}$ and $e_2^\text{MC}$ are, respectively, the sampling and discretization error, in the mean squared sense. For the root mean squared error (RMSE) to be lower than a prescribed tolerance $\epsilon$, it is sufficient to limit both $e_1^\text{MC}$ and $e_2^\text{MC}$ to, at most, $\epsilon^2/2$. 

For $\hat{F}_{h,M}^\text{MC}$ and all other CDF estimators considered in this work, we assume that the number of interpolation points $S+1$ is large enough for the interpolation error to be negligible and for $e$ to represent the error between the estimator and the true CDF $F$ of $Q$.

\subsection{Standard Multilevel Monte Carlo (MLMC) without smoothing}
\label{subsec:MLMC}

Rather than sampling $Q_M$ on a single spatial mesh, one considers a sequence of approximations $Q_{M_l}$ ($l=0,\dots,L_{\text{max}}$) of $Q$ associated with discrete meshes \{$\mathcal{T}_{M_l}, l=0,\dots,L_{\text{max}}$\}. Here $M_l$ the number of cells in mesh $\mathcal{T}_{M_l}$ and $M_{l-1}=2^{-d}M_l$, where $d$ is the spatial dimension. The idea behind this approach is to start by performing cheap-to-compute samples on a coarse mesh, and then gradually correct the resulting estimate of $F_{h,M}$ by sampling on finer grids, where generating a realization is more computationally expensive. One rewrites $\mathbb{E}[\mathcal{I}_n(Q_M)]$ as a telescopic sum
\begin{subequations}
\begin{align}
  \mathbb{E}[\mathcal{I}_n(Q_M)] = \mathbb{E}[\mathcal{I}_n(Q_0)] + \sum_{l=0}^{L_\text{max}}\mathbb{E}[\mathcal{I}_n(Q_{M_l}) - \mathcal{I}_n(Q_{M_{l-1}})]\equiv \sum_{l=0}^{L_\text{max}} \mathbb{E}[\mathcal{I}_n(Y_l)],  
\end{align}
where $\mathcal{I}_n(Y_l)$ ($l=0,\dots,L_\text{max}$) are given by
\begin{align}\label{I_Y_l}
  \mathcal{I}_n(Y_l) = 
  \begin{cases}   
     \mathcal{I}_n(Q_{M_l}) - \mathcal{I}_n(Q_{M_{l-1}}) & 1\leq l\leq L_{\text{max}} \\
     \mathcal{I}_n(Q_{M_l}) & l=0.
  \end{cases}
\end{align} 
\end{subequations}
The MLMC estimator for $\mathbb{E}[\mathcal{I}_n(Q_M)]$ is defined as
\begin{align}\label{MLMC_ind}
   \hat{\mathcal{I}}_{n,M}^\text{MLMC} = \sum_{l=0}^{L_{\text{max}}} \hat{\mathcal{I}}_n^\text{MC}(Y_l) =  
  \sum_{l=0}^{L_{\text{max}}} \frac{1}{N_l} \sum_{j=1}^{N_l} \mathcal{I}_n(Y_l^{(j)}).
\end{align} 
While $\mathbb{V}[\mathcal{I}_n(Q_{M_l})]$ remains approximately constant with $l$, $\mathbb{V}[\mathcal{I}_n(Y_l)]$ decreases with $l$, allowing the estimator $\hat{\mathcal{I}}_{n,M}^\text{MLMC}$ to have the same overall sampling error as its MC counterpart $\hat{\mathcal{I}}_{n,M}^\text{MC}$ using a decreasing number of samples $N_l$ as $l$ increases.

The MSE of the MLMC estimator $\hat{F}_{h,M}^\text{MLMC}$ for $F_{h,M}$ is bounded by
\begin{align}\label{MSE_MLMC}
  \mathbb{E}[\norm{e(\hat{F}_{h,M}^\text{MLMC})}_{\infty}^2] &\leq \underbrace{\mathbb{E}[\norm{\hat{F}_{h,M}^\text{MLMC} - \mathbb{E}[\hat{F}_{h,M}^\text{MLMC}]}_{\infty}^2]}_{e_1^\text{ML}} + \underbrace{\norm{F_{h,M} - F_h}^2}_{e_2^\text{ML}} \\\nonumber 
  &\leq \max_{0\leq n\leq S} \sum_{l=0}^{L_{\text{max}}} N_l^{-1} \mathbb{V}[\mathcal{I}_n(Y_l)] + \max_{0\leq n\leq S} |\mathbb{E}[\mathcal{I}_n(Q_{M_{L_{\text{max}}}}) - \mathcal{I}_n(Q)]|^2,
\end{align}              
where $e_1^\text{ML}$ and $e_2^\text{ML}$ are, respectively, the sampling and discretization error, in the mean squared sense. From the triangle inequality it follows that
\begin{align}\label{discret_cond}
  \max_{0\leq n\leq S} |\mathbb{E}[\mathcal{I}_n(Y_{L_{\text{max}}})| \approx \max_{0\leq n\leq S} |\mathbb{E}[\mathcal{I}_n(Q_{M_{L_{\text{max}}}}) - \mathcal{I}_n(Q)]|.
\end{align}
Hence, to determine the maximum level $L_\text{max}$ of an MLMC simulation with given tolerance $\epsilon$, we check if $\max_{0\leq n\leq S} |\mathbb{E}[\mathcal{I}_n(Y_{L})|\leq \epsilon/\sqrt{2}$ is satisfied for the current level $L$. Once the value of $L_\text{max}$ is found, we can compare the performance of MLMC to MC by performing the latter on this finest level, i.e., for $M$ in Section~\ref{subsec:MC} equal to $M_{L_{\text{max}}}$, re-using the samples already computed with MLMC at this level. Since $M=M_{L_{\text{max}}}$, this strategy ensures that the bound for $e_2^\text{ML}$ in~\eqref{MSE_MLMC} is the same as the bound for $e_2^\text{MC}$ in~\eqref{MSE_MC}.
To achieve an RSME error of at most $\epsilon$, it is sufficient that $e_1^\text{ML}\leq \epsilon^2/2$ and $e_2^\text{ML}\leq \epsilon^2/2$. 

\subsection{Standard Multilevel Monte Carlo (MLMC) with smoothing}
\label{subsec:MLMC_smooth}

The jump discontinuity in the indicator function may lead to a slow decay of $\mathbb{V}[\mathcal{I}_n(Y_l)]$ and make MLMC slower than MC for sufficiently large values of the error tolerance $\epsilon$~\cite{Lu2016}. To accelerate the variance decay and to improve the computational efficiency of MLMC, a sigmoid-type smoothing function $g$ can be used to remove the singularity in the indicator function. We consider two different smoothing functions, namely a polynomial proposed by Giles et al.~\cite{Giles2015} (which we will refer to as $g_\text{G}$) and the CDF of a Gaussian kernel in the context of kernel density estimation~\cite{Rosenblatt1956} (denoted by $g_\text{K}$). 

\subsubsection{Smoothing via Giles' polynomial}
\label{subsec:smooth_Giles}

Under the assumption, made in Section~\ref{sec:MLMC}, that the PDF $f_M$ of $Q_M$ is at least $d$-times continuously differentiable on $[a-\xi_0,b+\xi_0]$ for some $d\in\mathbb{N}_0$ and $\xi_0>0$, we define a smoothing function $g_\text{G}:\mathbb{R}\rightarrow \mathbb{R}$ that satisfies~\cite{Giles2015}
\begin{enumerate}
 \item cost of computing $g_\text{G}(s)\leq C$ $\forall s\in\mathbb{R}$ for some constant $C$
 \item $g_\text{G}$ is Lipschitz continuous
 \item $g_\text{G}(s)=1$ for $s<-1$ and $g_\text{G}(s)=0$ for $s>1$
 \item $\int_{-1}^1 s^k \mathcal{I}_{]-\infty,0]}(s - g_\text{G}(s)) \text ds=0$ for $k=0,\dots,d-1$.
\end{enumerate}
The function $g_\text{G}$ can be constructed as the uniquely determined polynomial of degree (at most) $d+1$, such that $\int_{-1}^1 s^k g_\text{G}(s) \text ds = (-1)^k / (k+1)$ for $k=0,\dots,d-1$, $g_\text{G}(1)=0$ and $g_\text{G}(-1)=1$. To obtain an appropriate smoothing function, $g_\text{G}$ is extended with
\begin{align}
  g_\text{G}(s) = 
  \begin{cases}   
     1 & s<-1 \\
     0 & s>1. 
  \end{cases}
\end{align} 

The indicator function $\mathcal{I}_n(Q_{M_l})$ at each level $l$ ($l=0,\dots,L_\text{max}$) is replaced by $g_\text{G}((Q_{M_l} - q_n)/\delta_{\text{G},l})$, where the \emph{bandwidth} $\delta_{\text{G},l}$ is a measure of the width over which the discontinuity in $\mathcal{I}_n(Q_{M_l})$ is smoothed out. The MLMC estimator with smoothing for $\mathbb{E}[\mathcal{I}_n(Q_M)]$ is given by
\begin{align}
   \hat{g}_{n,M}^\text{MLMC} = \sum_{l=0}^{L_{\text{max}}} \hat{g}_n^\text{MC}(Y_l), 
\end{align}
where $\hat{g}_n^\text{MC}(Y_l) = N_l^{-1} \sum_{j=1}^{N_l} g_n(Y_l^{(j)})$ with 
\begin{align}
  g_n(Y_l^{(j)}) = 
  \begin{cases}   
     g_\text{G}\left(\frac{Q_{M_l}^{(j)} - q_n}{\delta_{\text{G},l}}\right) - g_\text{G}\left(\frac{Q_{M_{l-1}}^{(j)} - q_n}{\delta_{\text{G},l}}\right) & 1\leq l\leq L_{\text{max}} \\
      g_\text{G}\left(\frac{Q_{M_l}^{(j)} - q_n}{\delta_{\text{G},l}}\right) & l=0.
  \end{cases}
\end{align} 
The MLMC estimator with smoothing $\hat{F}_{h,\delta_\text{G},M}^\text{MLMC}$ for $F_{h,M}$ is given by a piecewise polynomial interpolation with degree $\max(d,1)$
\begin{align}\label{MLMC_CDF_Giles}
  \hat{F}_{h,\delta_\text{G},M}^\text{MLMC}(q) = \sum_{n=0}^S \hat{g}_{n,M}^\text{MLMC}\ \phi_n(q) = \sum_{n=0}^S \sum_{l=0}^{L_{\text{max}}} \hat{g}_n^\text{MC}(Y_l) \phi_n(q).
\end{align}  
The MSE of $\hat{F}_{h,\delta_\text{G},M}^\text{MLMC}$ is bounded by
\begin{align}\label{MSE_MLMC_smooth}
   \mathbb{E}[\norm{e(\hat{F}_{h,\delta_\text{G},M}^\text{MLMC})}_{\infty}^2]\leq &\underbrace{\mathbb{E}[\norm{\hat{F}_{h,\delta_\text{G},M}^\text{MLMC} - \mathbb{E}[\hat{F}_{h,\delta_\text{G},M}^\text{MLMC}]}_{\infty}^2]}_{e_1^\text{ML}} + \underbrace{\norm{F_{h,M} - F_h}^2}_{e_2^\text{ML}} \\\nonumber
   &+\underbrace{\norm{\mathbb{E}[\hat{F}_{h,\delta_\text{G},M}^\text{MLMC}] - F_{h,M}}_{\infty}^2}_{e_3^\text{ML}}.
\end{align}
Compared to~\eqref{MSE_MLMC},~\eqref{MSE_MLMC_smooth} contains an additional term $e_3^\text{ML}$, which is the (mean squared) smoothing error. To achieve an RSME error of at most $\epsilon$, it is sufficient that $e_1^\text{ML}\leq \epsilon^2/4$, $e_2^\text{ML}\leq \epsilon^2/2$ and $e_3^\text{ML}\leq \epsilon^2/4$. 

Finding the optimal value for the bandwidth $\delta_{\text{G},l}$ at each level $l$ ($l=0,\dots,L_\text{max}$) is crucial to the performance of MLMC. The error $e_3^\text{ML}$ should be as close as possible to $\epsilon^2/4$ in order to maximize the smoothing and, hence, the variance decay with increasing level. A possible strategy consists of the following steps.
\begin{enumerate}
   \item Given the error tolerance $\epsilon$, at level $l=0$ estimate $\delta_{\text{G},l,n}$ for each interpolation point $q_n$ in $S_h=\{q_n,n=0,\dots,S\}$ by solving
    \begin{align}
      \frac{1}{N_l^0}\left|\sum_{j=1}^{N_l^0} \left[g_\text{G}\left(\frac{Q_{M_l}^{(j)} - q_n}{\delta_{\text{G},l,n}} \right) - \mathcal{I}_n(Q_{M_l}^{(j)}) \right]\right| = \frac{\epsilon}{2}
    \end{align}
    based on a set of initial samples $\{Q_{M_l}^{(j)}\}_{j=1}^{N_l^0}$.
    \item Define the smoothing parameter for level $l=0$, $\delta_{\text{G},l}$, as
          \begin{align}
             \delta_{\text{G},l} = \max_{0\leq n\leq S} \delta_{\text{G},l,n}.
          \end{align}
    \item Repeat steps 1 and 2 for each new level.      
\end{enumerate}        

We follow the numerical algorithm in~\ref{sec:Implementation} to compute $\hat{F}_{h,\delta_\text{G},M}^\text{MLMC}$ and to measure the associated computational cost (see Section~\ref{subsec:cost}). This algorithm is inspired by the approaches~\cite{Cliffe2011, Lu2016}. To compare the performance of the MLMC simulation with that of the corresponding MC simulation at the highest level $L_\text{max}$ (which ensures that $e_2^\text{MC}\leq \epsilon^2/2$), we also compute the number of MC samples of $Q_{L_\text{max}}$ required to satisfy $e_1^\text{MC}\leq \epsilon^2/2$ and the resulting computational cost of the MC estimator $\hat{F}_{h,M}^\text{MC}$ (see~\ref{subsec:MLMC_imp}).

\subsubsection{Smoothing based on kernel density estimation (KDE)}
\label{subsec:smooth_KDE}

We propose an alternative way of smoothing the indicator function. It is grounded in Kernel Density Estimation (KDE), a nonparametric estimation of the PDF of a random variable. Let $q_1,\dots,q_n$ be independent and identically distributed samples drawn from the distribution of a QoI $Q$ with unknown PDF $f$. Then the KDE of $f$ is given by
\begin{align}
   \hat{f}_{\delta_\text{K}}(q) = \frac{1}{n}\sum_{i=1}^n \mathcal K_{\delta_\text{K} }(q-q_i) = \frac{1}{n\delta_\text{K}}\sum_{i=1}^n \mathcal K\left(\frac{q-q_i}{\delta_\text{K}}\right)
\end{align}
where $\mathcal K$ is a \emph{kernel} and $\delta_\text{K}$ is the \emph{bandwidth}. We refer to $K_{\delta_\text{K}}(z)=(1/\delta_\text{K}) \mathcal K(z/\delta_\text{K})$ for $z \in \mathbb{R}$ as a \emph{scaled kernel}. For a Gaussian kernel, the KDE of $f$ is 
\begin{align}
   \hat{f}_{\delta_\text{K}}(q) = \frac{1}{n\delta_\text{K}} \sum_{i=1}^n \frac{1}{\sqrt{2\pi}} \exp\left[-\frac{1}{2}\left(\frac{q-q_i}{\delta_\text{K}}\right)^2 \right].
\end{align}
A corresponding estimate of the CDF $F$ is then obtained by considering the CDF of the Gaussian kernel, which yields
\begin{align}\label{KDE_CDF}
   \hat{F}_{\delta_\text{K}}(q) = \frac{1}{n}\sum_{i=1}^n \Phi\left(\frac{q-q_i}{\delta_\text{K}}\right).
\end{align}
Here $\Phi$ is the CDF of the standard normal distribution and plays the role of indicator smoothing function. The bandwidth $\delta_\text{K}$ is the counterpart of $\delta_\text{G}$ defined in Section~\ref{subsec:smooth_Giles}. Based on~\eqref{KDE_CDF}, at each level $l$ ($l=0,\dots,L_\text{max}$) we replace $\mathcal{I}_n(Q_{M_l})$ with $\Phi[(q_n - Q_{M_l})/\delta_{\text{K},l}] \equiv g_\text{K}[ (q_n - Q_{M_l})/\delta_{\text{K},l}]$ and define an MLMC estimator with smoothing for $F_{h,M}$ similar to the one given by~\eqref{MLMC_CDF_Giles} but with the smoothing function $g_\text{K}$. The methodology of Section~\ref{subsec:smooth_Giles} is then used to bound its MSE and find the optimal value of the bandwidth $\delta_{\text{K},l}$ at each level $l$ ($l=0,\dots,L_\text{max}$). The algorithm in~\ref{subsec:MLMC_imp} is deployed to compute $\hat{F}_{h,\delta_\text{K},M}^\text{MLMC}$ and to measure its  computational cost.

\subsection{Stratified Multilevel Monte Carlo (sMLMC)}
\label{subsec:sMLMC}

In stratified Monte Carlo (sMC), one divides the domain $\mathcal{D}$ of the uncertain input parameter $W$ into $r$ mutually exclusive and exhaustive regions $\mathcal{D}_i$ ($i=1,\dots,r$) called \emph{strata}. Let $p_i=\mathbb{P}(W\in \mathcal{D}_i)$ be the probability of $W$ being in stratum $\mathcal{D}_i$. The expectation value of $Q_M(W)$ is
\begin{align}
  \mathbb{E}[Q_M(W)]=\sum_{i=1}^r p_i \zeta_i,
\end{align}
where $\zeta_i$ ($i=1,\dots,r$) is the $i^{\text{th}}$ stratum mean given by
\begin{align}
  \zeta_i = \int_{\mathcal{D}_i} Q_M(w)\ \mathrm d F_{W}^{(i)}(w) = \frac{1}{p_i} \int_{\mathcal{D}_i} Q_M(w)\ \mathrm d F_W(w).
\end{align}
Here $F_{W}^{(i)}$ is the conditional CDF of $W$ given that $W\in\mathcal{D}_i$. The sMC estimator for $\mathbb{E}[Q_M(W)]$ is given by
\begin{align}\label{sMC_exp}
   \hat{Q}_M^\text{sMC} = \sum_{i=1}^r \frac{p_i}{n_i}\sum_{j=1}^{n_i} Q_M^{(j,i)},
\end{align} 
where $n_i$ is the number of independent samples generated in the $i$th stratum for each $i=1,\dots,r$ with $\sum_{i=1}^{r} n_i\equiv N$, and $Q_M^{(j,i)}$ is the $j$th sample of $Q_M$ that has a corresponding input parameter ($W$) in stratum $i$. The variance of this estimator is 
\begin{align}\label{sMC_var}
  \mathbb{V}[\hat{Q}_M^\text{sMC}] = \sum_{i=1}^r \frac{\sigma_i^2 p_i^2}{n_i}, \qquad \sigma_i^2 = \int_{\mathcal{D}_i} (Q_M(w) - \zeta_i)^2 \mathrm d F_W^{(i)}(w)
\end{align}
with $\sigma_i^2$ being the variance of $Q_M$ within stratum $i$. 

Two common choices for $n_i$ are proportional and optimal allocations~\cite{Fishman1996}. For proportional allocation, $n_i=Np_i$ and 
\begin{align}
   \hat{Q}_M^\text{sMC, p} &= \frac{1}{N}\sum_{i=1}^r \sum_{j=1}^{n_i} Q^{(j,i)}\\
   \mathbb{V}[\hat{Q}_M^\text{sMC, p}] &= \frac{1}{N}\sum_{i=1}^r \sigma_i^2 p_i = \mathbb{V}[\hat{Q}_M^\text{MC}] - \frac{1}{N}\sum_{i=1}^r p_i (\zeta_i - \mathbb{E}(Q_M))^2\label{prop_alloc},
\end{align} 
 which shows that stratification produces an estimator with a lower variance than its MC counterpart. For optimal allocation, $n_i=\alpha_i N$ with $\alpha_i = \sigma_i p_i/(\sum_{k=1}^r \sigma_k p_k)$ for $i=1,\dots,r$, yielding
\begin{align}
  \mathbb{V}[\hat{Q}_M^\text{sMC, o}] &= \frac{1}{N}\left(\sum_{i=1}^r \sigma_i p_i\right)^2.
\end{align}
It can be shown~\cite{Fishman1996} that this is the smallest variance possible for an sMC estimator, which, given~\eqref{prop_alloc}, is also smaller than the variance of the corresponding MC estimator.

Replacing $Q_M$ with $\mathcal{I}_n(Q_M)$ in~\eqref{sMC_exp} leads to
\begin{align}\label{sMC_cdf}
   \hat{\mathcal{I}}_{n,M}^\text{sMC} = \sum_{i=1}^r \frac{p_i}{n_i}\sum_{j=1}^{n_i} \mathcal{I}_n(Q_M^{(j,i)}).
\end{align}
To combine the benefits of stratification with those of the multigrid approach, we replace MC at each level of the MLMC algorithm with sMC and refer to the resulting algorithm as ``stratified Multilevel Monte Carlo'' (sMLMC). In analogy to~\eqref{MLMC_ind}, the sMLMC estimator for $\mathbb{E}[\mathcal{I}_n(Q_M)]$ is defined as
\begin{align}\label{sMLMC_ind}
   \hat{\mathcal{I}}_{n,M}^\text{sMLMC} &= \sum_{l=0}^{L_{\text{max}}} \hat{\mathcal{I}}_n^\text{sMC}(Y_l) \\\nonumber
   &= \sum_{l=0}^{L_{\text{max}}} \sum_{i=1}^r \frac{p_i}{n_{i,l}}\sum_{j=1}^{n_i} \mathcal{I}_n(Y_l^{(j,i)}),
\end{align} 
where $\mathcal{I}_n(Y_l)$ with $l=0,\dots,L_\text{max}$ are defined in~\eqref{I_Y_l}. The MSE of the sMLMC estimator $\hat{F}_{h,M}^\text{sMLMC}$ for $F_{h,M}$ is bounded by
\begin{align}\label{MSE_sMLMC}
  \mathbb{E}[\norm{e(\hat{F}_{h,M}^\text{sMLMC})}_{\infty}^2] &\leq \underbrace{\mathbb{E}[\norm{\hat{F}_{h,M}^\text{sMLMC} - \mathbb{E}[\hat{F}_{h,M}^\text{sMLMC}]}_{\infty}^2]}_{e_1^\text{sML}} + \underbrace{\norm{F_{h,M} - F_h}^2}_{e_2^\text{sML}} \\\nonumber 
  &\leq \max_{0\leq n\leq S} \sum_{l=0}^{L_{\text{max}}} \mathbb{V}[\hat{\mathcal{I}}_n^\text{sMC}(Y_l)] + \max_{0\leq n\leq S} |\mathbb{E}[\mathcal{I}_n(Q_{M_{L_{\text{max}}}}) - \mathcal{I}_n(Q)]|^2
\end{align}        
where the bound for $e_2^\text{sML}$ is the same as the bound for $e_2^\text{MC}$ in~\eqref{MSE_MC}, provided $M_{L_{\text{max}}}=M$. 
To reduce the computational cost further, we again smooth the indicator function, yielding an additional error term in~\eqref{MSE_sMLMC} similar to $e_3^\text{ML}$ in~\eqref{MSE_MLMC_smooth}. The resulting sMLMC estimator with smoothing for $\mathbb{E}[\mathcal{I}_n(Q_M)]$ is defined as
\begin{align}
   \hat{g}_{n,M}^\text{sMLMC} = \sum_{l=0}^{L_{\text{max}}} \hat{g}_n^\text{sMC}(Y_l), 
\end{align}
where $\hat{g}_n^\text{sMC}(Y_l) = \sum_{i=1}^r (p_i / n_{i,l}) \sum_{j=1}^{n_{i,l}} g_n(Y_l^{(j,i)})$ and 
\begin{align}
  g_n(Y_l^{(j,i)}) = 
  \begin{cases}   
     g_\text{G}\left(\frac{Q_{M_l}^{(j,i)} - q_n}{\delta_{\text{G},l}}\right) - g_\text{G}\left(\frac{Q_{M_{l-1}}^{(j,i)} - q_n}{\delta_{\text{G},l}}\right) & 1\leq l\leq L_{\text{max}} \\
     g_\text{G}\left(\frac{Q_{M_l}^{(j,i)} - q_n}{\delta_{\text{G},l}}\right) & l=0
  \end{cases}
\end{align} 
for polynomial-based smoothing, or
\begin{align}
  g_n(Y_l^{(j,i)}) = 
  \begin{cases}   
     g_\text{K}\left(\frac{q_n - Q_{M_l}^{(j,i)}}{\delta_{\text{K},l}}\right) - g_\text{K}\left(\frac{q_n - Q_{M_{l-1}}^{(j,i)}}{\delta_{\text{K},l}}\right) & 1\leq l\leq L_{\text{max}} \\
     g_\text{K}\left(\frac{q_n - Q_{M_l}^{(j,i)}}{\delta_{\text{K},l}}\right) & l=0
  \end{cases}
\end{align} 
for KDE-based smoothing.
The sMLMC estimator with smoothing for $F_{h,M}$ is then given by
\begin{align}
  \hat{F}_{h,\delta,M}^\text{sMLMC}(q) = \sum_{n=0}^S \hat{g}_{n,M}^\text{sMLMC}\ \phi_n(q) = \sum_{n=0}^S \sum_{l=0}^{L_{\text{max}}} \hat{g}_n^\text{sMC}(Y_l) \phi_n(q).
\end{align} 
with $\delta$ given by $\delta_\text{G}$ or $\delta_\text{K}$. To compute $\hat{F}_{h,\delta_\text{G},M}^\text{sMLMC}$ and $\hat{F}_{h,\delta_\text{K},M}^\text{sMLMC}$ and measure the associated computational cost, we deploy the algorithm in~\ref{subsec:sMLMC_imp}.

\subsection{Relative cost of MLMC and sMLMC versus MC}
\label{subsec:cost}

We estimate the total cost of computing the MLMC estimator without smoothing of $F_{h,M}$ as an average over $N_\text{real}$ independent realizations of the corresponding algorithm,
\begin{align}
\mathcal{C}(\hat{F}_{h,M}^\text{MLMC})=\frac{1}{N_\text{real}}\sum_{k=1}^{N_\text{real}}\sum_{l=0}^{L_{\text{max},k}^\text{MLMC}} \bar{w}_{l}^{(k)} N_{l}^{(k)},
\end{align}
where $\bar{w}_{l}^{(k)}$ is the average computational cost of computing a sample of $Q_{M_l}$ on level $l$ for realization $k$, and $L_{\text{max},k}^\text{MLMC}$ denotes the finest discretization level at which the sampling is performed for this realization. 

For sMLMC without smoothing,
\begin{align}
  \mathcal{C}(\hat{F}_{h,M}^\text{sMLMC}) = \frac{1}{N_\text{real}} \sum_{k=1}^{N_\text{real}} \sum_{l=0}^{L_{\text{max},k}^\text{sMLMC}}\sum_{i=1}^r \bar{w}_{i,l}^{(k)} n_{i,l}^{(k)},
\end{align}
where $w_{i,l}^{(k)}$ is the average computational cost for computing a sample of $Q_{M_l}$ in stratum $i$ on level $l$ for realization $k$, and $L_{\text{max},k}^\text{sMLMC}$ refers to the finest discretization level for this realization. 

To compare the performance of MLMC and sMLMC with or without smoothing with that of MC, we consider $N_\text{real}$ realizations of the MC algorithm and perform the $k^{\text{th}}$ realization on level $L_{\text{max},k}^\text{MLMC}$ of the corresponding realization of MLMC without smoothing. This provides a single average cost,
\begin{align}
  \mathcal{C}(\hat{F}_{h,M}^\text{MC}) = \frac{1}{N_\text{real}}\sum_{k=1}^{N_\text{real}}\bar{w}_{L_{\text{max},k}^\text{MLMC}} N_\text{MC}^{(k)}
\end{align}
against which to compare $\mathcal{C}(\hat{F}_{h,M}^\text{MLMC})$, $\mathcal{C}(\hat{F}_{h,\delta,M}^\text{MLMC})$, $\mathcal{C}(\hat{F}_{h,M}^\text{sMLMC})$ and $\mathcal{C}(\hat{F}_{h,\delta,M}^\text{sMLMC})$ (with $\delta=\delta_\text{G}$ for polynomial-based smoothing or $\delta=\delta_\text{K}$ for KDE-based smoothing). Here $N_\text{MC}^{(k)}$ is the number of samples computed in the $k^{\text{th}}$ realization of the MC algorithm.

\section{Numerical results}
\label{sec:Results}

We consider two testbed problems: linear diffusion with a random diffusion coefficient, and inviscid Burgers' equation with a random initial condition. In both cases, the random input variable $W$ is drawn from a truncated lognormal PDF $f_W(w;\mu,\sigma,w_l,w_u)$ defined on $[w_l,w_u]$ with mean $\mu$ and variance $\sigma^2$, i.e., 
\begin{align}
   f_W = \frac{\sqrt{2}}{\sqrt{\pi}\sigma w}
   \begin{cases}   
      \dfrac{\exp(-\frac{(\ln w - \mu)^2}{2\sigma^2})}{ \text{erf}\left(\frac{\ln w_u - \mu}{\sqrt{2}\sigma}\right) - \text{erf}\left(\frac{\ln w_l - \mu}{\sqrt{2}\sigma}\right)} & \text{for } w\in[w_l,w_u] \\
      0 & \text{otherwise}.
   \end{cases}
\end{align}
We assume the PDF $f$ of the QoI $Q(W)$ to be at least 3 times continuously differentiable such that we can use a cubic-spline interpolation to approximate its CDF $F$. To compare the performance of MC, MLMC and sMLMC, we measure their average computational cost over $N_\text{real}=50$ independent runs for error tolerances $\epsilon=0.01$, 0.008 and 0.005. The highest possible maximum level $L_\text{max}^{\star}$ for MLMC/sMLMC we consider is $L_\text{max}^{\star}=7$. At each level we use proportional allocation for the sMC algorithm. We choose an identical number of warmup samples $N_l^0$ ($l=0,\dots,L_\text{max}^{\star}$) for all independent runs of the MLMC algorithm, but use different values for MLMC with and without smoothing to eliminate as much as possible any oversampling ($N_l^0$ is lower when smoothing is applied). For sMLMC, we define the number of warmup samples $n_{i,l}^0$ in each stratum $i=1,\dots,r$, based on $N_l^0$ with $l=0,\dots,L_\text{max}^{\star}$ and the chosen allocation strategy (identical for all independent sMLMC runs). The values of $N_l^0$ for sMLMC is typically lower than those for MLMC and different between the cases with and without smoothing, again to avoid oversampling.

For each value of $\epsilon$, the $k^{\text{th}}$ run involves the following steps.
\begin{enumerate}
  \item Perform non-smoothed MLMC, yielding a maximum level $L_{\text{max},k}^\text{MLMC}$.
  \item Perform MC with $M=M_{L_{\text{max},k}^\text{MLMC}}$, re-using already computed samples from the MLMC run.
  \item Perform smoothed MLMC using Giles' polynomial with a computed smoothing parameter $\delta_{\text{G},l}$ at each level $l$.
  \item Perform smoothed MLMC using KDE with a computed smoothing parameter $\delta_{\text{K},l}$ at each level $l$.
  \item Perform non-smoothed sMLMC.
  \item Perform smoothed sMLMC using KDE.
\end{enumerate}

\subsection{Linear diffusion equation}
\label{subsec:LD}

Consider 
\begin{subequations}\label{diff_eqn}
  \begin{align}
      &\frac{\partial u}{\partial t} = D\frac{\partial^2 u}{\partial x^2}, \quad x\in(0,4), \quad t>0  \\
      &u(0,t)=-1, \quad u(4,t)=1 \\
      &u(x,0)=\tanh \left(\frac{2-x}{0.05}\right)
   \end{align}\label{diff}    
\end{subequations}
Here $D\in\Omega_D=[1,4]$ is a lognormal random variable with PDF $f_D(w;3,3,1,4)$. Our goal is to compute the CDF of a QoI 
\begin{align}
   Q=10\int_0^4 \mathrm dx\ u^2(x,0.2).        
\end{align}  
We discretize~\eqref{diff} in space using a central finite difference scheme, and then apply the Crank-Nicholson method to the resulting system of initial value problems. The matrix associated with this linear system is tridiagonal Toeplitz, hence we apply the Thomas algorithm to solve it at reduced computational complexity. This numerical scheme is second-order accurate in both space and time. We approximate the CDF $F$ of Q on the interval [14,28], which contains the support of $F$, and use $S+1=29$ interpolation points, i.e., set $h=0.5$.

\begin{figure}[tphb]
\begin{center}
\includegraphics[width=0.49\textwidth]{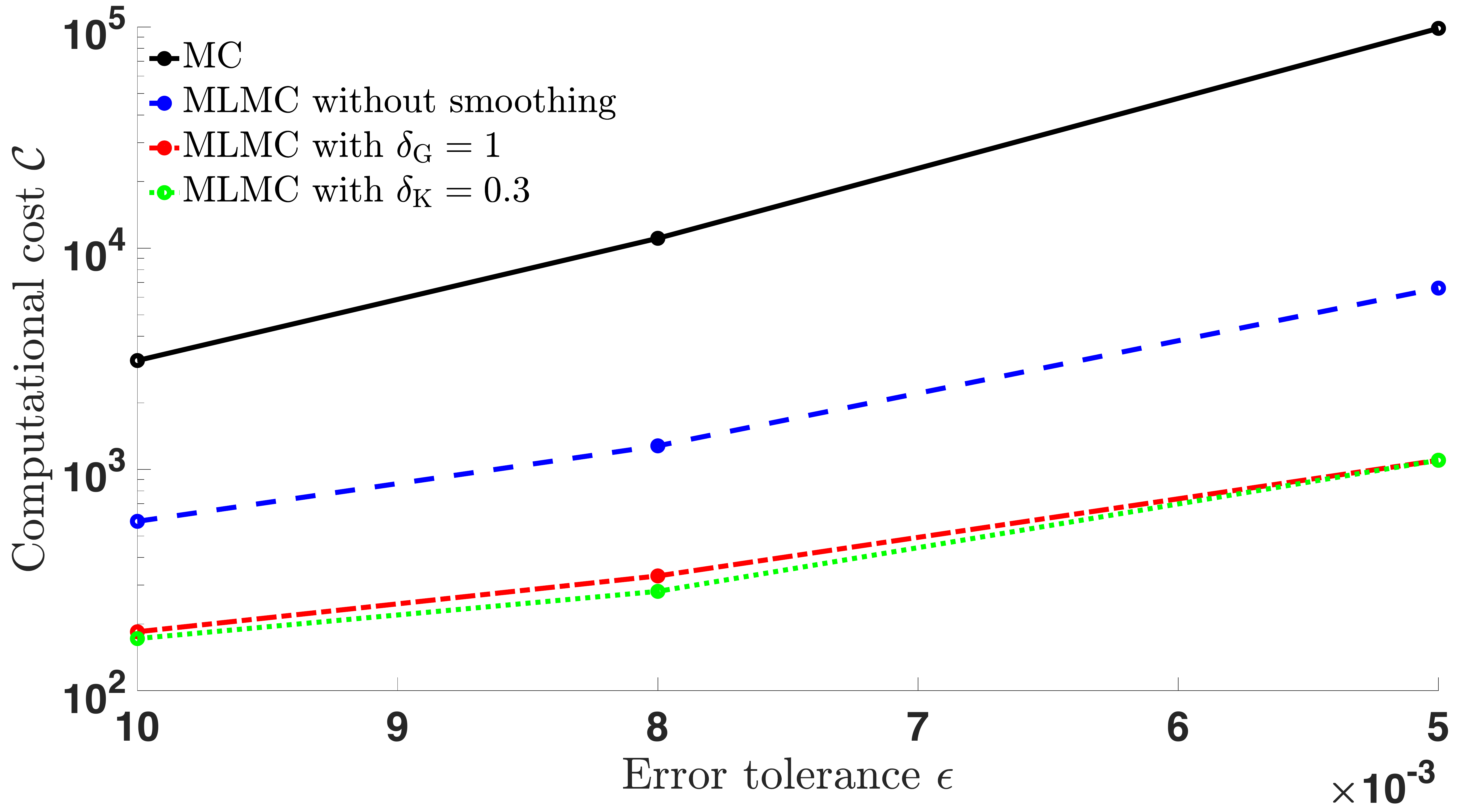}
\includegraphics[width=0.49\textwidth]{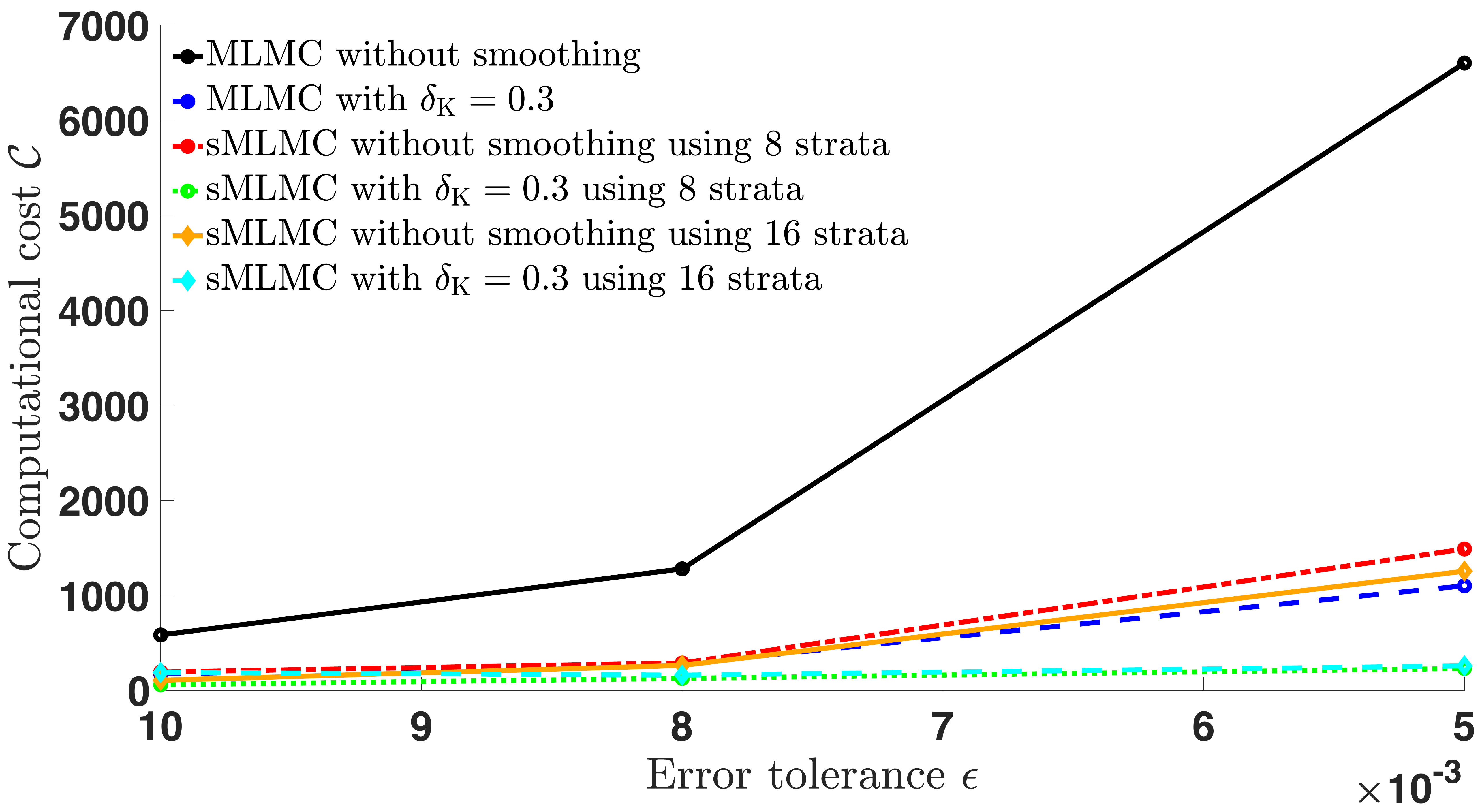}
\end{center}
\caption{Computational cost of MC, MLMC without smoothing and with polynomial-based or KDE-based smoothing (left), and of sMLMC without smoothing and with KDE-based smoothing (right) using 8 or 16 strata for the linear diffusion problem.}
\label{MC_MLMC_comp_LD}
\end{figure}

We first compare the computational cost of MC, $\mathcal{C}(\hat{F}_{h,M}^\text{MC})$, to that of non-smoothed MLMC, $\mathcal{C}(\hat{F}_{h,M}^\text{MLMC})$, and MLMC with the polynomial-based, $\mathcal{C}(\hat{F}_{h,\delta_\text{G},M}^\text{MLMC})$, and KDE-based, $\mathcal{C}(\hat{F}_{h,\delta_\text{K},M}^\text{MLMC})$, smoothing of the indicator function. Figure~\ref{MC_MLMC_comp_LD} (left) shows that even without smoothing the cost of MLMC is about an order of magnitude lower than that of MC, at the lowest error tolerance considered ($\epsilon=0.005$). Applying the smoothing further reduces the computational cost by almost an order of magnitude at this tolerance (Fig.~\ref{MC_MLMC_comp_LD}, right). We find that KDE-based smoothing yields a lower computational cost than its polynomial-based counterpart, and hence only consider KDE for smoothing the indicator function in the sMLMC algorithm.

\begin{figure}[tphb]
\begin{center}
\includegraphics[width=0.49\textwidth]{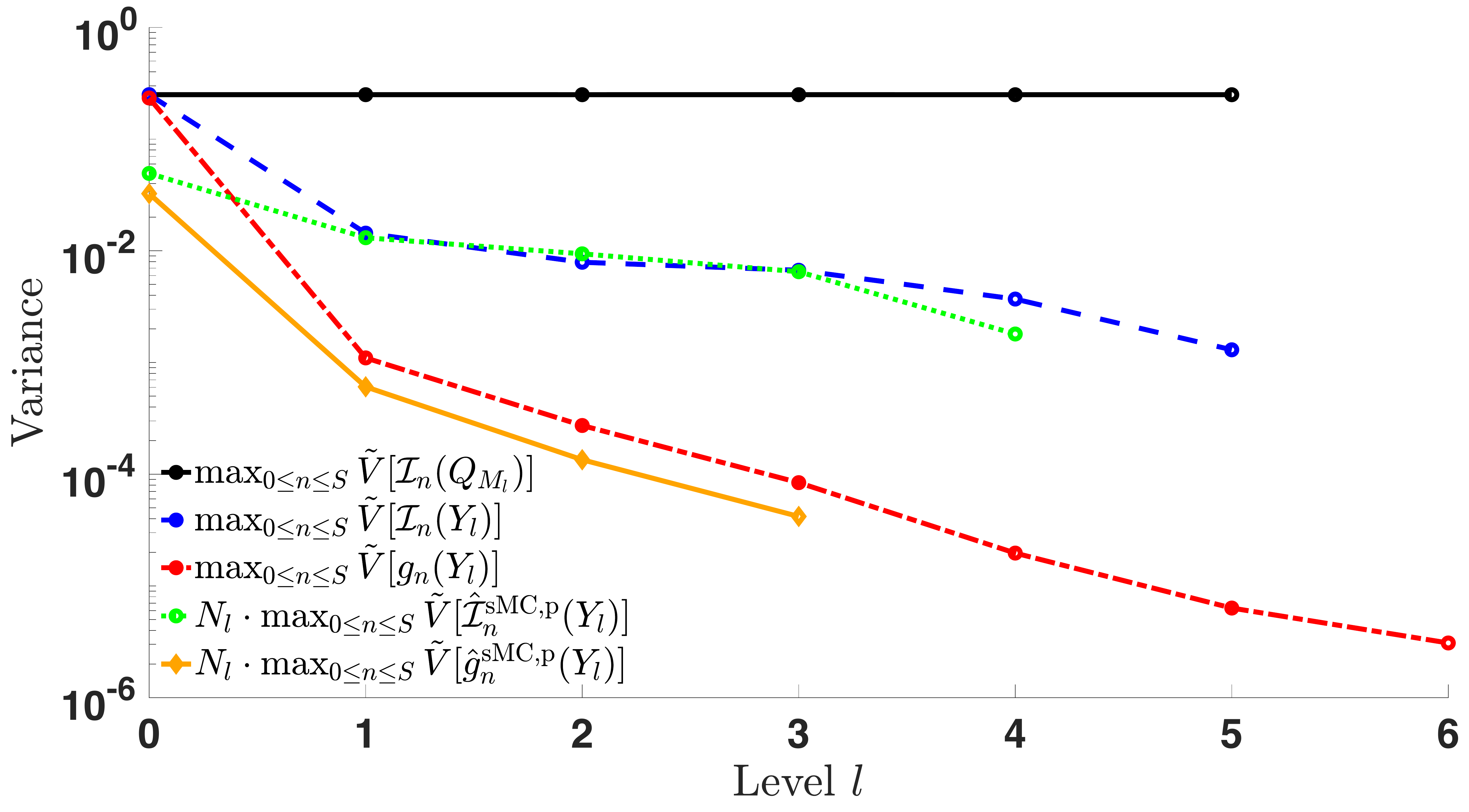}
\includegraphics[width=0.49\textwidth]{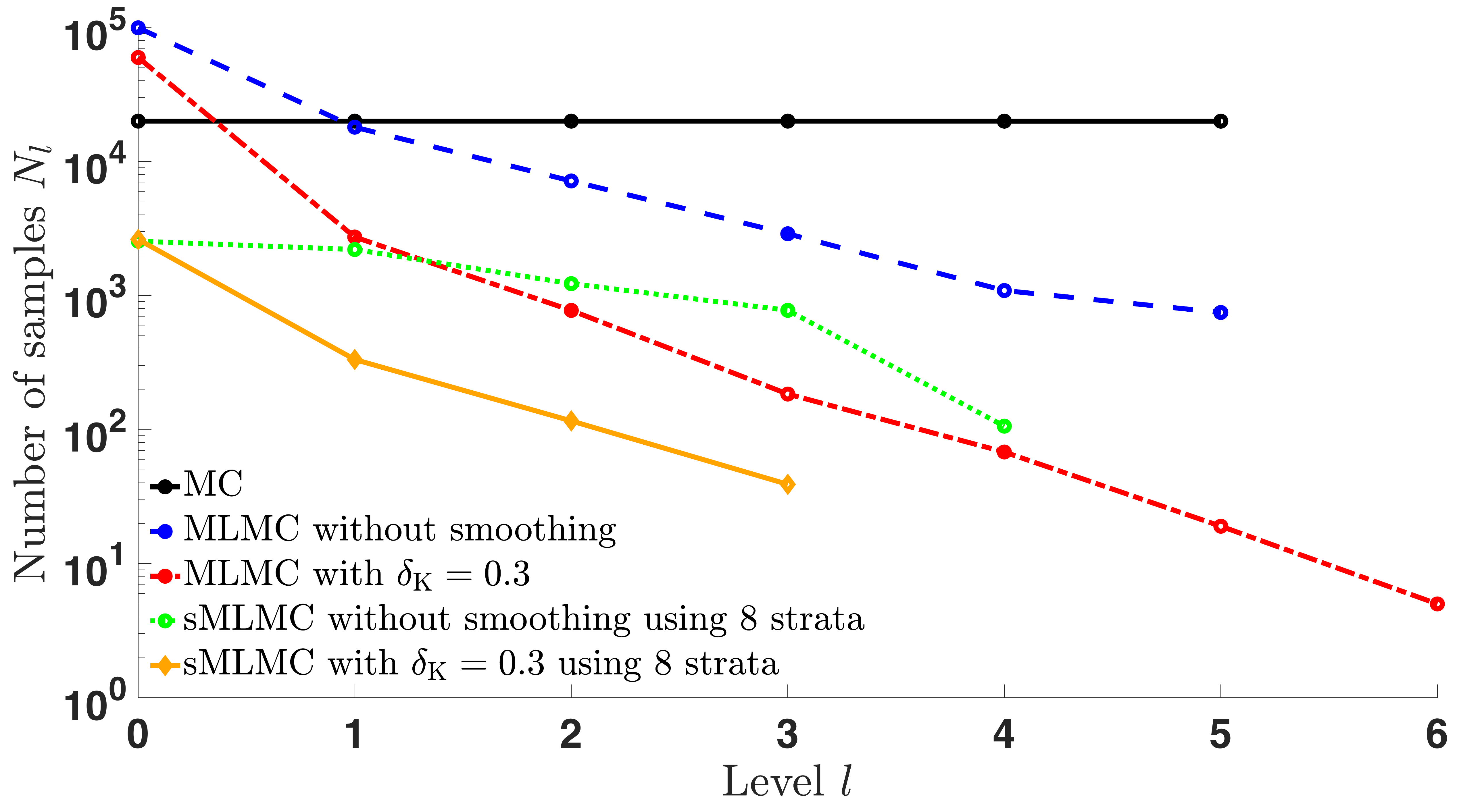}
\end{center}
\caption{Evolution of the variance (a) and number of samples (b) with level for a single run of the linear diffusion problem and $\epsilon=0.005$. For the sMLMC runs we use 8 strata.}
\label{MC_MLMC_comp_LD2}
\end{figure}

Figure~\ref{MC_MLMC_comp_LD2} (left) illustrates the evolution of $\tilde{V}[\mathcal{I}_n(Q_{M_l})]$, $\tilde{V}[\mathcal{I}_n(Y_l)]$ and $\tilde{V}[g_n(Y_l)]$ (using KDE-based smoothing) with level $l$ for a single run and $\epsilon=0.005$. Here $\tilde{V}$ denotes a sample estimate of $\mathbb{V}$. While $\tilde{V}[\mathcal{I}_n(Q_{M_l})]$ remains approximately constant as the spatial resolution increases, $\tilde{V}[\mathcal{I}_n(Y_l)]$ and $\tilde{V}[g_n(Y_l)]$ decay as the spatial mesh is refined, so that fewer samples are needed at higher levels of discretization (Fig.~\ref{MC_MLMC_comp_LD2}, right).

Next, we stratify the sample space of $D$, $\Omega_D$, into 8 or 16 strata of equal width and compare the computational cost of the (non-stratified) MLMC algorithm without smoothing, $\mathcal{C}(\hat{F}_{h,M}^\text{MLMC})$, to that of sMLMC without smoothing, $\mathcal{C}(\hat{F}_{h,M}^\text{sMLMC})$ and with KDE-based smoothing, $\mathcal{C}(\hat{F}_{h,\delta_\text{K},M}^\text{sMLMC})$, for 8 and 16 strata and the error tolerances considered.  Figure~\ref{MC_MLMC_comp_LD} (right) demonstrates that stratifying the non-smoothed MLMC algorithm yields similar computational cost savings as applying smoothing to the indicator function. Combining KDE-based smoothing with stratification further reduces the computational cost by almost an order of magnitude for $\epsilon=0.005$, yielding between one and two orders of magnitude cost savings compared to the non-smoothed MLMC method at this tolerance. We note that for 16 strata and $\epsilon=0.01$, smoothing increases, rather than decreases, the overall computational cost of the sMLMC algorithm. In this case, the number of samples on each level is already quite low for non-smoothed sMLMC and smoothing the indicator function is not beneficial. As predicted by~\eqref{prop_alloc}, Figure~\ref{MC_MLMC_comp_LD2} (left) demonstrates that stratification yields an additional variance reduction. Here and for the smoothed case, we compare the quantity $N_l\cdot\max_{0\leq n\leq S}\tilde{V}[\hat{\mathcal{I}}_n^\text{sMC, p}(Y_l)]$ with $\tilde{V}[\mathcal{I}_n(Y_l)]$. Combining stratification with smoothing yields further variance reduction. This is translated into a decrease in the required number of samples (Fig.~\ref{MC_MLMC_comp_LD2}, right).

\subsection{Inviscid Burgers' equation}
\label{subsec:IB}

Consider 
\begin{subequations}
\begin{align}
  &\frac{\partial u}{\partial t} + \frac{1}{2}\frac{\partial u^2}{\partial x}=0, \quad x\in(0,2), \quad t>0\label{Burgers_inviscid_eqn} \\
  &u(0,t)=2, \quad  u(2,t)=0 \\
  &u(x,0)=
  \begin{cases}
     U_1 & 0<x\leq 1 \\
     0 & 1<x<2.
  \end{cases}\label{U_1}
\end{align}\label{Burgers_inviscid}
\end{subequations}
The initial state $U_1\in\Omega_{U_1}=[0,2]$ is a lognormal random variable with PDF $f_{U_1}(w;1.5,1,0,2)$. Our goal is to compute the CDF of a QoI  
\begin{align}
   Q=10\int_0^2 \mathrm dx\ u^2(x,0.5).
\end{align}  
We discretize~\eqref{Burgers_inviscid} using the Godunov method, which is first-order accurate in both space and time. This is a conservative finite volume scheme which solves a Riemann problem at each inter-cell boundary. We approximate the CDF $F$ of Q over the interval [15,65], which contains the support of $F$, and use $S+1=101$ interpolation points, i.e., set $h=0.5$.

\begin{figure}[tphb]
\begin{center}
\includegraphics[width=0.49\textwidth]{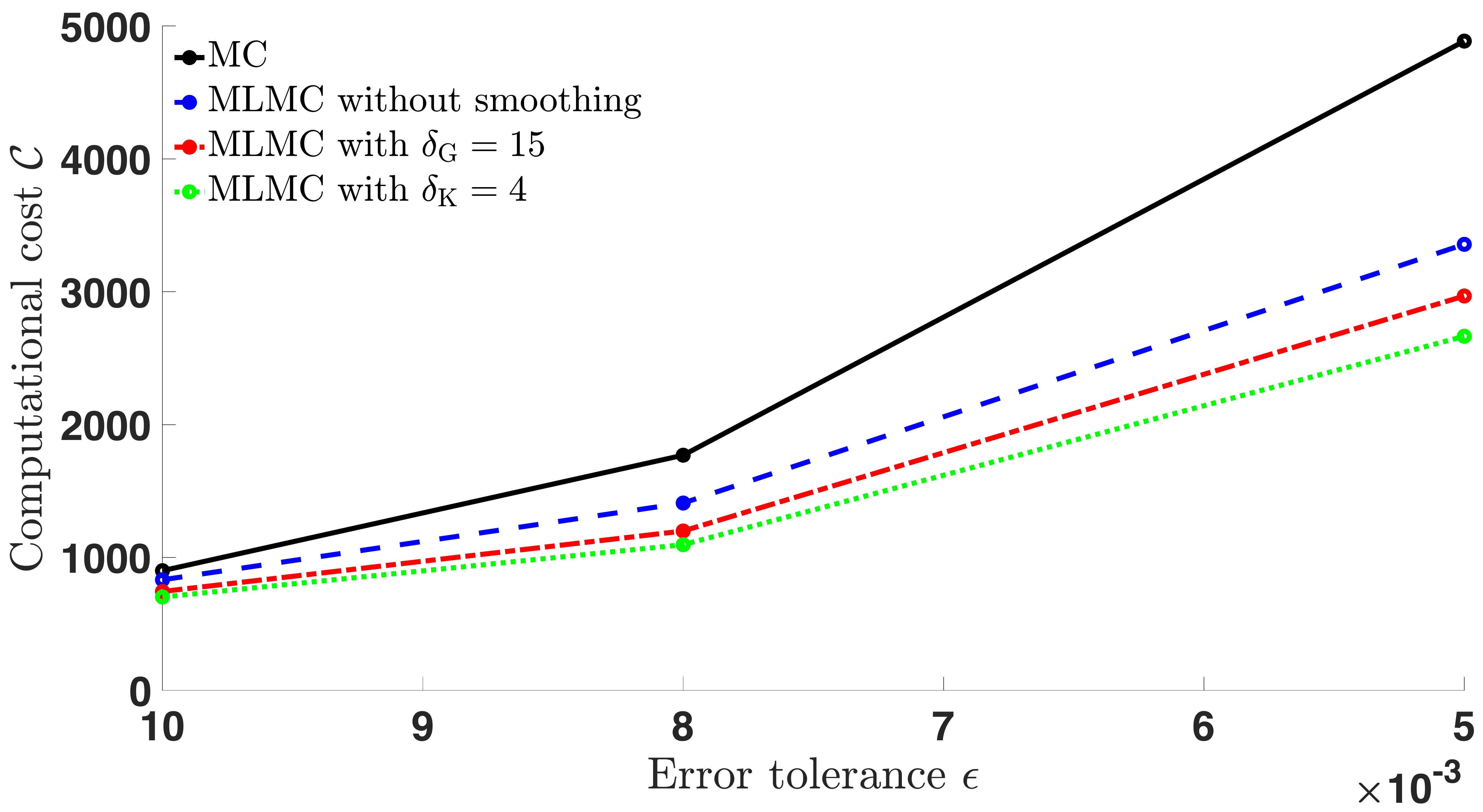}
\includegraphics[width=0.49\textwidth]{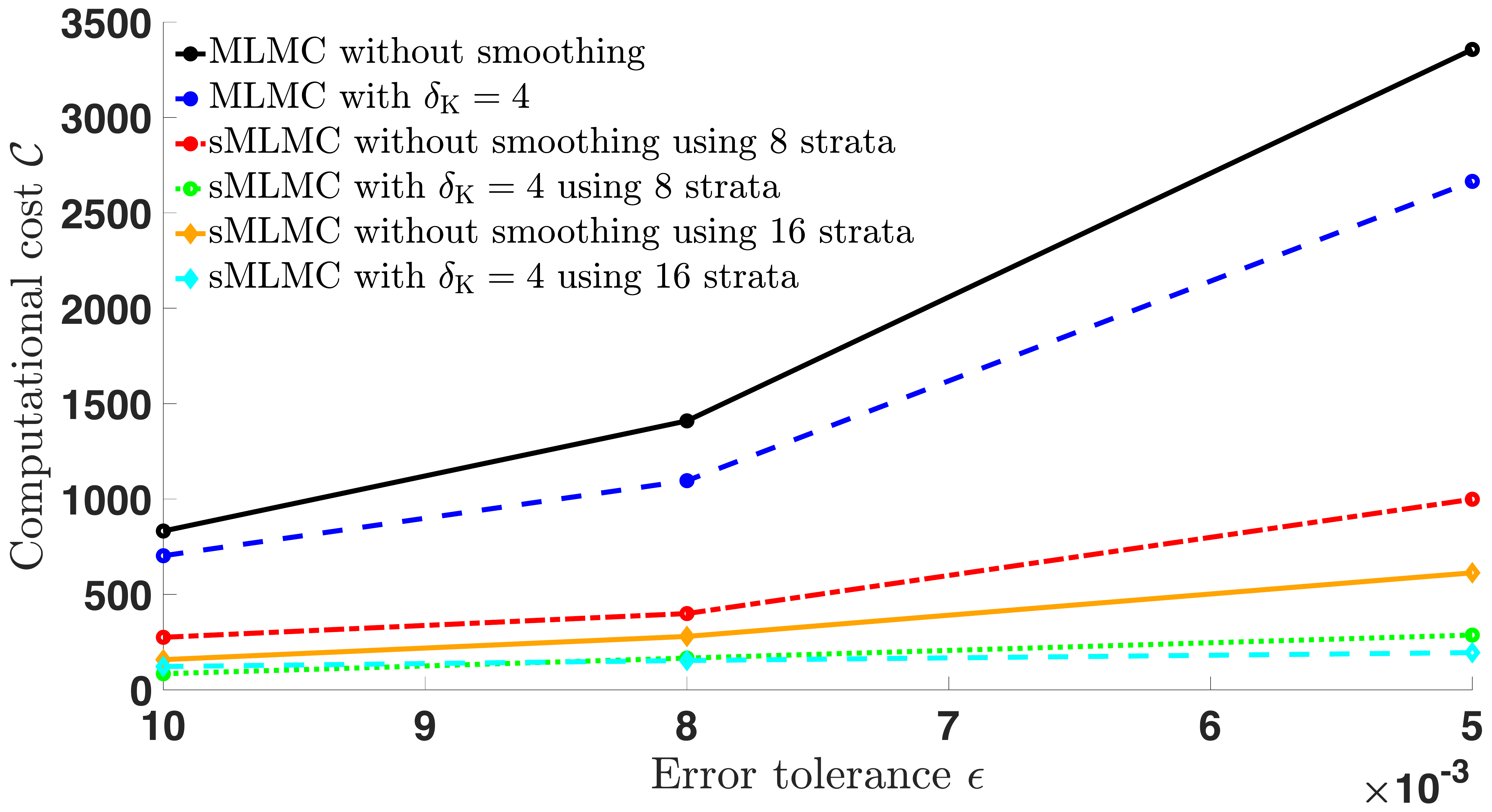}
\end{center}
\caption{Computational cost of MC, MLMC without smoothing and with polynomial-based or KDE-based smoothing (left), and of MLMC without smoothing/with KDE-based smoothing, and sMLMC without smoothing/with KDE-based smoothing (right), using 8 or 16 strata for the inviscid Burgers' problem.}
\label{MC_MLMC_comp_IB}
\end{figure}

We compare the computational cost of MC, $\mathcal{C}(\hat{F}_{h,M}^\text{MC})$, to that of non-smoothed MLMC, $\mathcal{C}(\hat{F}_{h,M}^\text{MLMC})$, and MLMC with polynomial-based, $\mathcal{C}(\hat{F}_{h,\delta_\text{G},M}^\text{MLMC})$, and KDE-based, $\mathcal{C}(\hat{F}_{h,\delta_\text{K},M}^\text{MLMC})$, smoothing of the indicator function. Figure~\ref{MC_MLMC_comp_IB} (left) illustrates that, as in the linear diffusion problem, MLMC is more efficient than MC even without smoothing. However, the computational cost savings are much lower in this case. The KDE-based smoothing is more efficient than the polynomial-based smoothing, especially at the lowest error tolerance $\epsilon=0.005$. However, the gain in computational efficiency from smoothing is relatively modest. 

Next, we stratify the sample space of $U_1$, $\Omega_{U_1}$, into 8 or 16 strata of equal width and compare the computational cost of the (non-stratified) MLMC algorithm without smoothing, $\mathcal{C}(\hat{F}_{h,M}^\text{MLMC})$, to that of sMLMC without smoothing, $\mathcal{C}(\hat{F}_{h,M}^\text{sMLMC})$, and with KDE-based smoothing, $\mathcal{C}(\hat{F}_{h,\delta_\text{K},M}^\text{sMLMC})$, for 8 and 16 strata and the error tolerances considered. Figure~\ref{MC_MLMC_comp_IB} (right) shows that the computational saving from applying stratification are much more substantial than those originating from smoothing the indicator function, unlike for the linear diffusion problem where the non-smoothed sMLMC and smoothed MLMC algorithm had a similar computational cost. Combining KDE-based smoothing with stratification again yields the most efficient result, saving about an order of magnitude in computational cost compared to standard MLMC with or without smoothing.

\section{Conclusions }
\label{sec:Concl}

We constructed a stratified Multilevel Monte Carlo (sMLMC) algorithm for estimating the cumulative distribution function (CDF) of a quantify of interest (QoI) in a problem with random input parameters or a random initial state. Our method combines the benefits of multigrid from standard MLMC with variance reduction from stratified sampling in each of the discretization levels. We also explored the use of Gaussian kernel density estimator (KDE) in lieu of the currently used polynomial-based smothers.

Our study yields the following major conclusions:
\begin{enumerate}
  \item For all test problems and error tolerances considered, the computational cost of non-smoothed sMLMC is smaller than that of non-smoothed MLMC. This can be attributed to the fact that, as spatial resolution increases, the variance decays faster for sMLMC than for MLMC due to the additional variance reduction from the stratification at each level.
  \item Stratifying MLMC can yield either similar or substantially higher computational cost savings compared to smoothing the indicator function, depending on the considered problem.
  \item Smoothing the indicator function in the sMLMC algorithm  yields noticeable additional computational cost savings, with a few exceptions at high tolerances where smoothing either produces a negligible improvement or even a reduction in the algorithm's efficiency. Like in the non-smoothed case, smoothed sMLMC is more efficient than its smoothed MLMC counterpart.
  \item For the same level of smoothing error, KDE-based smoothing yields a more efficient algorithm than its polynomial-based counterpart. The gain in efficiency depends on the required error tolerance and the problem considered.
\end{enumerate}

\section{Acknowledgments}

This work was supported in part by Defense Advanced Research Project Agency under award number 101513612, by Air Force Office of Scientific Research under award number FA9550-17-1-0417, and by U.S. Department of Energy under award number DE-SC0019130.

\appendix

\section{Numerical algorithms for computing $\hat{F}_{h,\delta,M}^\text{MLMC}$ and $\hat{F}_{h,\delta,M}^\text{sMLMC}$}
\label{sec:Implementation}

\subsection{MLMC with smoothing}
\label{subsec:MLMC_imp}

\begin{algorithm}[H]
{\scriptsize
\caption{Standard Multilevel Monte Carlo with smoothing}
\SetAlgoLined
\SetKwData{Left}{left}\SetKwData{This}{this}\SetKwData{Up}{up}
\SetKwFunction{Union}{Union}\SetKwFunction{FindCompress}{FindCompress}
\SetKwInOut{Input}{Input}\SetKwInOut{Output}{Output}\SetKwInOut{Procedure}{Procedure}

\Input{The RMSE accuracy $\epsilon$, $S+1$ interpolation points $a=q_0 < q_1 < \dots < q_S=b$ with separation $h=(b-a)/S$; a sequence of discrete meshes \{$\mathcal{T}_{M_l}, l=0,\dots,L_{\text{max}}$\}; and the initial number of samples $N_l^0$ at each level $l$\;}
\Output{An estimate of the CDF $F(q)$\;}
\Procedure: 
Initialize $L=-1$\;
\While{$L< L_{\mathrm{max}}$}{
       Set $L=L+1$\;
       Draw $N_L^0$ samples of the random input parameter (IP)/initial condition (IC) ($\star$)\; 
       \eIf{$L=0$}{
            Compute $N_0^0$ samples of $Q_{M_0}$ based on ($\star$)\;
       }{     
            Compute $N_L^0$ samples of $Q_{M_L}$ and $Q_{M_{L-1}}$ based on ($\star$)\;
       }
       Compute $\delta_L$\; 
       \For{$n=0,\dots,S$}{
            \For{$j=1,\dots,N_L^0$}{
                 Compute $g_n(Y_L^{(j)})$\;
            }    
       }
       Compute the computational cost at level $L$, $\bar{w}_L$\;
       \For{$n=0,\dots,S$}{
            Compute $\hat{\mathcal{I}}_n^\text{MC}(Y_L)$ and $\hat{g}_n^\text{MC}(Y_L)$\;       
            Compute $\tilde{V}[g_n(Y_L)]=\sum_{j=1}^{N_L^0} (g_n(Y_L^{(j)}) - \hat{g}_n^\text{MC}(Y_L))^2 / N_L^0$ 
       }
       See next page
}
}
\end{algorithm} 

\begin{algorithm}
{\scriptsize
\SetAlgoLined
Set $N_L=\displaystyle\text{ceil}\left(\max_{0\leq n\leq S} 4\epsilon^{-2}\sqrt{\tilde{V}[g_n(Y_L)]/\bar{w}_L}\left(\sum_{k=0}^L \sqrt{\tilde{V}[g_n(Y_k)]\bar{w}_k}\right)\right)$\; 
Draw $\max(N_L - N_L^0,0)$ samples of the random IP/IC ($\dag$)\; 
       \eIf{$L=0$}{
            Compute $\max(N_0 - N_0^0,0)$ samples of $Q_{M_0}$ based on ($\dag$)\;
       }{     
            Compute $\max(N_L - N_L^0,0)$ samples of $Q_{M_L}$ and $Q_{M_{L-1}}$ based on ($\dag$)\;
       }  
\For{$n=0,\dots,S$}{
     \For{$j=N_L^0+1,\dots,N_L$}{
          Compute $g_n(Y_L^{(j)})$\;
     }
}          
Compute the computational cost at level $L$, $\bar{w}_L$\;
\For{$n=0,\dots,S$}{
     Compute $\hat{\mathcal{I}}_n^\text{MC}(Y_L)$ and $\hat{g}_n^\text{MC}(Y_L)$\;       
     Compute $\tilde{V}[g_n(Y_L)]=\sum_{j=1}^{N_L} (g_n(Y_L^{(j)}) - \hat{g}_n^\text{MC}(Y_L))^2 / N_L$\;
}
Set $N_L^{\star}=N_L$\;
\SetAlgoLined
\For{$l=0,\dots,L-1$}{
     Set $N_l=\displaystyle\text{ceil}\left(\max_{0\leq n\leq S} 4\epsilon^{-2}\sqrt{\tilde{V}[g_n(Y_l)]/\bar{w}_l}\left(\sum_{k=0}^L \sqrt{\tilde{V}[g_n(Y_k)]\bar{w}_k}\right)\right)$\;
     Draw $\max(N_0 - N_0^{\star},0)$ samples of the random IP/IC ($\ddag$)\; 
     \eIf{$l=0$}{
          Compute $\max(N_0 - N_0^{\star},0)$ samples of $Q_{M_0}$ based on ($\ddag$)\;
     }{
          Compute $\max(N_l - N_l^{\star},0)$ samples of $Q_{M_l}$ and $Q_{M_{l-1}}$ based on ($\ddag$)\;
     }     
     \For{$n=0,\dots,S$}{
          \For{$j=N_l^{\star}+1,\dots,N_l$}{
               Compute $g_n(Y_l^{(j)})$\;
          }
     }      
     Compute the computational cost at level $l$, $\bar{w}_l$\;
     \For{$n=0,\dots,S$}{
          Compute $\hat{\mathcal{I}}_n^\text{MC}(Y_l)$ and $\hat{g}_n^\text{MC}(Y_l)$\;       
          Compute $\tilde{V}[g_n(Y_l)]=\sum_{j=1}^{N_l} (g_n(Y_l^{(j)}) - \hat{g}_n^\text{MC}(Y_l))^2 / N_l$\;
     }
     Set $N_l^{\star}=N_l$\;
}  
\If{($L\geq 1$ $\mathrm{and}$ $\displaystyle\max_{0\leq n\leq S}|\hat{\mathcal{I}}_n^\mathrm{MC}(Y_L)|\leq \epsilon/\sqrt{2}$) $\mathrm{or}$ ($L=L_{\mathrm{max}}$) }{
     Compute the cost of MLMC, $\mathcal{C}(\hat{F}_{h,\delta,M}^\text{MLMC})$\;
     Compute the MLMC estimator of $F(q)$, $\hat{F}_{h,\delta,M}^\text{MLMC}(q)$\;
     Set $N_\text{MC}=2\epsilon^{-2}\max_{0\leq n\leq S} \tilde{V}[\mathcal{I}_n(Q_{M_{L_\text{max}}})]$\;
     Compute the cost of MC, $\mathcal{C}(\hat{F}_{h,M}^\text{MC})$\;
     Compute $\max(N_\text{MC} - N_{L_\text{max}},0)$ samples of $Q_{M_{L_\text{max}}}$\;
     Compute the MC estimator of $F(q)$, $\hat{F}_{h,M}^\text{MC}(q)$\;
} 
}
\end{algorithm}

\subsection{sMLMC with smoothing}
\label{subsec:sMLMC_imp}

\begin{algorithm}[H]
{\scriptsize
\caption{Stratified MLMC with smoothing}
\SetAlgoLined
\SetKwData{Left}{left}\SetKwData{This}{this}\SetKwData{Up}{up}
\SetKwFunction{Union}{Union}\SetKwFunction{FindCompress}{FindCompress}
\SetKwInOut{Input}{Input}\SetKwInOut{Output}{Output}\SetKwInOut{Procedure}{Procedure}

\Input{Input parameters of Algorithm in Section~\ref{subsec:MLMC_imp}; the strata $\mathcal{D}_i$ with $i=1,\dots,r$; and their probabilities $p_i$\;}
\Output{An estimate of the CDF $F(q)$\;}
\Procedure
\\Initialize $L=-1$\;
 \While{$L< L_{\mathrm{max}}$}{
  Set $L=L+1$\;
  For $i=1,\dots,r$, define initial number of samples in stratum $i$ at level $L$, $n_{i,L}^0$, based on $N_L^0$ and chosen allocation strategy, and draw $n_{i,L}^0$ samples of random IP/IC from stratum $i$ ($\star$)\; 
  \eIf{$L=0$}{
       \For{$i=1,\dots,r$}{
            Compute $n_{i,0}^0$ samples of $Q_{M_0}$ based on ($\star$)\;
       }
  }{        
       \For{$i=1,\dots,r$}{
            Compute $n_{i,L}^0$ samples of $Q_{M_L}$ and $Q_{M_{L-1}}$ based on ($\star$)\; 
       } 
  }      
  Compute $\delta_L$ with the combined samples from all strata\;
  \For{$i=1,\dots,r$}{
       \For{$j=1,\dots,n_{i,L}^0$}{
            \For{$n=0,\dots,S$}{
                 Compute $g_n(Y_L^{(j,i)})$\;
            }
       }
       Compute sample variance estimate $\tilde{V}[g_n(Y_L^{(i)})]$ of $\mathbb{V}[g_n(Y_L^{(i)})]$\;
  }   
  \For{$i=1,\dots,r$}{   
       Compute the average computational cost per sample in stratum $i$ at level $L$, $\bar{w}_{i,L}$\;
  }    
  See next page
  }
  \SetAlgoLined
Set $n_{i,L}^{\star}=n_{i,L}$ for all $i=1,\dots,r$\;
\For{$n=0,\dots,S$}{
     \For{$i=1,\dots,r$}{
          Compute $n_{i,L,n}=4\epsilon^{-2}\sqrt{\tilde{V}[g_n(Y_L^{(i)})] p_i^2 / \bar{w}_{i,L}} \displaystyle\sum_{k=0}^L\displaystyle\sum_{i=1}^r\sqrt{\tilde{V}[g_n(Y_k^{(i)})] p_i^2 \bar{w}_{i,k}}$\;
     }
}   
\For{$i=1,\dots,r$}{
     Set $n_{i,L}=\displaystyle\text{ceil}\left(\max_{0\leq n\leq S} n_{i,L,n}\right)$\;
}
Set $N_L=\displaystyle\sum_{i=1}^r n_{i,L}$\;
See next page
}
\end{algorithm} 

\SetNlSty{texttt}{(}{)}
\begin{algorithm}
{\scriptsize
Draw $n_{i,L} - n_{i,L}^0$ samples of random IP/IC from stratum $i$\;
\eIf{$L=0$}{
     \For{$i=1,\dots,r$}{
          Draw $n_{i,0} - n_{i,0}^0$ samples of $Q_{M_0}$ based on ($\dag$)\;
     }
}{        
     \For{$i=1,\dots,r$}{
          Draw $n_{i,L} - n_{i,L}^0$ samples of $Q_{M_L}$ and $Q_{M_{L-1}}$ based on ($\dag$)\;
     } 
}     
\For{$i=1,\dots,r$}{
     \For{$j=n_{i,L}^0+1,\dots,n_{i,L}$}{
          \For{$n=0,\dots,S$}{
               Compute $g_n(Y_L^{(j,i)})$\;
          }
     }
     Compute improved estimate $\tilde{V}[g_n(Y_L^{(i)})]$ of $\mathbb{V}[g_n(Y_L^{(i)})]$\;
}         
\SetAlgoLined
\For{$i=1,\dots,r$}{   
     Compute the average computational cost per sample in stratum $i$ at level $L$, $\bar{w}_{i,L}$\;
} 
\For{$n=0,\dots,S$}{
     Compute $\hat{\mathcal{I}}_n^\text{sMC}(Y_L)$ and $\hat{g}_n^\text{sMC}(Y_L)$\;       
}
Set $n_{i,L}^{\star}=n_{i,L}$ for all $i=1,\dots,r$\;
\For{$l=0,\dots,L-1$}{
\For{$n=0,\dots,S$}{
     \For{$i=1,\dots,r$}{
          Compute $n_{i,l,n}=4\epsilon^{-2}\sqrt{\tilde{V}[g_n(Y_l^{(i)})] p_i^2 / \bar{w}_{i,l}} \displaystyle\sum_{k=0}^L\displaystyle\sum_{i=1}^r\sqrt{\tilde{V}[g_n(Y_k^{(i)})] p_i^2 \bar{w}_{i,k}}$\;
     }
}   
\For{$i=1,\dots,r$}{
     Set $n_{i,l}=\displaystyle\text{ceil}\left(\max_{0\leq n\leq S} n_{i,l,n}\right)$\;
}
Set $N_l=\displaystyle\sum_{i=1}^r n_{i,l}$\;   
Draw $n_{i,L} - n_{i,L}^{\star}$ samples of random IP/IC from stratum $i$\;
\eIf{$l=0$}{
     \For{$i=1,\dots,r$}{
          Compute $n_{i,0} - n_{i,0}^{\star}$ samples of $Q_{M_0}$ based on ($\ddag$)\;
     }
}{        
     \For{$i=1,\dots,r$}{
          Compute $n_{i,l} - n_{i,l}^{\star}$ samples of $Q_{M_l}$ and $Q_{M_{l-1}}$ based on ($\ddag$)\;
     } 
}     
See next page (part (a))
}  
See next page (part (b))
}
\end{algorithm}

\SetNlSty{texttt}{(}{)}
\begin{algorithm}
{\scriptsize
\SetAlgoLined
(a)\\
\For{$i=1,\dots,r$}{
     \For{$j=n_{i,l}^{\star}+1,\dots,n_{i,l}$}{
          \For{$n=0,\dots,S$}{
               Compute $g_n(Y_l^{(j,i)})$\;
          }
     }
     Compute improved estimate $\tilde{V}[g_n(Y_L^{(i)})]$ of $\mathbb{V}[g_n(Y_L^{(i)})]$\;
}      
\For{$i=1,\dots,r$}{   
     Compute the average computational cost per sample in stratum $i$ at level $l$, $\bar{w}_{i,l}$\;
}   
\For{$n=0,\dots,S$}{
     Compute $\hat{\mathcal{I}}_n^\text{sMC}(Y_l)$ and $\hat{g}_n^\text{sMC}(Y_l)$\;       
}
Set $n_{i,l}^{\star}=n_{i,l}$ for all $i=1,\dots,r$\;
\ (b)\\
\If{($L\geq 1$ $\mathrm{and}$ $\displaystyle\max_{0\leq n\leq S}|\hat{\mathcal{I}}_n^\mathrm{sMC}(Y_L)|\leq \epsilon/\sqrt{2}$) $\mathrm{or}$ ($L=L_{\mathrm{max}}$) }{
      Compute the cost of sMLMC, $\mathcal{C}(\hat{F}_{h,\delta,M}^\text{sMLMC})$\;
      Compute the sMLMC estimator of $F(q)$, $\hat{F}_{h,\delta,M}^\text{sMLMC}(q)$; 
} 
}
\end{algorithm}

\bibliography{paper}

\end{document}